\documentclass{article}
\catcode`\@=11
\newtheorem{lemma}{Lemma}
\newtheorem{proposition}{Proposition}

\newtheorem{example}{Example}
\newtheorem{remark}{Remark}

\def\demo{\noindent{\bf Proof .-}}
\def\section{\@startsection {section}{1}{\z@}{-3.5ex plus -1ex
minus-.2ex}{2.3ex plus .2ex}{\normalsize\bf}}

\begin{document}
\begin{center}
{\Large\bf \textsc{On the arithmetical rank of certain monomial ideals}}\footnote{MSC 2000: 13A15; 13F55, 14M10.}
\end{center}
\vskip.5truecm
\begin{center}
{Margherita Barile\footnote{Partially supported by the Italian Ministry of Education, University and Research.}\\ Dipartimento di Matematica, Universit\`{a} di Bari, Via E. Orabona 4,\\70125 Bari, Italy}\footnote{e-mail: barile@dm.uniba.it, Fax: 0039 080 596 3612}
\end{center}
\vskip1truecm
\noindent
{\bf Abstract} We determine a new technique which allows the computation of the arithmetical rank of certain monomial ideals.
\vskip0.5truecm
\noindent
Keywords: Arithmetical rank, monomial ideals, set-theoretic complete intersections.  

\section*{Introduction and Preliminaries} 
Given a commutative ring with identity $R$, the {\it arithmetical rank} of an ideal $I$ of $R$, denoted ara\,$I$, is defined as the minimum number of elements which generate $I$ {\it up to radical}, i.e., generate an ideal which has the same radical as $I$. Determining this number is, in general, a very hard open problem; a trivial lower bound is given by the height of $I$, but this is the actual value of ara\,$I$ only in special cases ($I$ is then called a {\it set-theoretic complete intersection}). There are, however, techniques which allow us to provide upper bounds. Some results in this direction have been recently proved in \cite{B5}, \cite{B6}, \cite{B7}, \cite{B4}, and less recently in \cite{B3}. These  apply especially to the case where $R$ is a polynomial ring over a field, and $I$ is a {\it monomial} ideal (i.e., an ideal generated by products of indeterminates) and are essentially based on the following criterion by Schmitt and Vogel
(see \cite{SV}, p.~249).
\begin{lemma}\label{lemma1} Let $P$ be a finite subset of elements of $R$. Let $P_0,\dots, P_r$ be subsets of $P$ such that
\begin{list}{}{}
\item[(i)] $\bigcup_{i=0}^rP_i=P$;
\item[(ii)] $P_0$ has exactly one element;
\item[(iii)] if $p$ and $p''$ are different elements of $P_l$ $(0<i\leq r)$ there is an integer $i'$ with $0\leq i'<i$ and an element $p'\in P_{i'}$ such that $pp''\in(p')$.
\end{list}
\noindent
Let $0\leq i\leq r$, and, for any $p\in P_i$, let $e(p)\geq1$ be an integer. We set $q_l=\sum_{p\in P_i}p^{e(p)}$. We will write $(P)$ for the ideal of $R$ generated by the elements of $P$.  Then we get
$$\sqrt{(P)}=\sqrt{(q_0,\dots,q_r)}.$$

\end{lemma}
If, in the construction given in the claim, we take all exponents $e(p)$ to be equal to 1, then $q_0,\dots, q_r$  are sums of generators; in this paper we present a new method, which gives rise to elements of the same form, but applies under a different assumption. It will enable us to determine the arithmetical rank of certain monomial ideals which could not be treated by the above lemma. \newline
For the determination of the arithmetical rank, every ideal can be replaced by its radical. In the sequel we will therefore throughout consider radical (or reduced) monomial ideals, i.e., ideals generated by squarefree monomials. These are the so-called {\it Stanley-Reisner ideals} of simplicial complexes.  For the basic notions on this topic we refer to \cite{BH}, Section 5. \newline
It is well-known that a reduced monomial ideal which is a set-theoretic complete intersection is Cohen-Macaulay; for the Stanley-Reisner ideals of one-dimensional simplicial complexes this condition is independent of the ground field, and is equivalent to the connectedness of the simplicial complex. It is not known whether all reduced monomial ideals which are Cohen-Macaulay are set-theoretic complete intersections. The question is open even for Gorenstein ideals. \newline
The problem of the arithmetical rank of monomial ideals has been intensively studied by several other authors over the past three decades: see \cite{G}, \cite{L1}, \cite{L2}, \cite{L3}, \cite{NV}, \cite{SchV}, \cite{SV}, \cite{Te} and \cite{Te2}. 
\section{The first result and some applications} Let $R$ be a commutative ring with identity. 
\begin{proposition}\label{main} Let $p_0, p_{11}, p_{12}, p_{21}, p_{22}\in R$ be such that $p_0$ divides $p_{11}p_{22}$, $p_{21}$ divides $p_{11}p_{12}$, and $p_{12}$ divides $p_{21}p_{22}$. Set $p_1=p_{11}+p_{12}$ and  $p_2=p_{21}+p_{22}$. Then 
$$\sqrt{(p_0, p_1, p_2)}=\sqrt{(p_0, p_{11}, p_{12}, p_{21}, p_{22})}.$$
\end{proposition}
\demo It suffices to show that $p_{11}, p_{12}, p_{21}, p_{22}\in\sqrt{(p_0, p_1, p_2)}.$ Let $a, b\in R$ be such that $p_{11}p_{22}=ap_0$ and $p_{11}p_{12}=bp_{21}$. Then
\begin{eqnarray*}
p_{11}^3&=&p_{11}^2(p_{11}+p_{12})-p_{11}^2p_{12}=p_{11}^2(p_{11}+p_{12})-p_{11}bp_{21}\\
&=&p_{11}^2(p_{11}+p_{12})-bp_{11}(p_{21}+p_{22})+bp_{11}p_{22}\\
&=&p_{11}^2(p_{11}+p_{12})-bp_{11}(p_{21}+p_{22})+bap_0\\
&=&p_{11}^2p_1-bp_{11}p_2+bap_0\in(p_0, p_1, p_2).
\end{eqnarray*}
\noindent
This shows that $p_{11}\in\sqrt{(p_0, p_1, p_2)}$, whence $p_{12}=p_1-p_{11}\in\sqrt{(p_0, p_1, p_2)}.$ The claim for $p_{21}$, $p_{22}$ follows by symmetry. This completes the proof. 
\par\medskip\noindent
The above result can be used  for computing the arithmetical rank of some monomial ideals. We first apply it for proving the set-theoretic complete intersection property of the Stanley-Reisner ideals of two simplicial complexes.  Recall that if the maximal faces of a simplicial complex on $N$ vertices all have the same cardinality $d$, then its Stanley-Reisner ideal has pure height $N-d$. Any Stanley-Reisner ideal that is Cohen-Macaulay has pure height (see \cite{BH}, Corollary 5.1.5), so that this holds, in particular, for set-theoretic complete intersections.    
\begin{example}\label{example1} {\rm In the polynomial ring $R=K[x_1,\dots, x_5]$, where $K$ is any field, consider the  ideal 
$$I=(x_1x_3,\ x_1x_4,\ x_2x_4,\ x_2x_5,\ x_3x_5),$$
\noindent
which is the Stanley-Reisner ideal of the Cohen-Macaulay simplicial complex on the vertex set $\{1,\dots, 5\}$ whose maximal faces are $\{1,2\}$,  $\{2,3\}$, $\{3,4\}$, $\{4,5\}$, $\{5,1\}$.
It is of pure height 3. Proposition \ref{main} applies to 
$$p_0=x_1x_3,\ p_{11}=x_1x_4,\ p_{12}=x_2x_5,\ p_{21}=x_2x_4,\ p_{22}=x_3x_5,$$
\noindent
so that
$$I=\sqrt{(x_1x_3,\ x_1x_4+x_2x_5,\ x_2x_4+x_3x_5)}.$$ 
\noindent
Hence ara\,$I=3$, and $I$ is a set-theoretic complete intersection. Note that three elements of $R$ generating $I$ up to radical cannot be found by applying Lemma \ref{lemma1} to the set of minimal monomial generators of $I$. }
\end{example}
\begin{example}{\rm In $R=K[x_1,\dots, x_6]$ consider the  ideal 
$$I=(x_1x_3,\ x_1x_4,\ x_1x_6,\ x_2x_4,\ x_2x_5,\ x_3x_5,\ x_3x_6,\ x_4x_6),$$
\noindent
which is the Stanley-Reisner ideal of the Cohen-Macaulay simplicial complex on the vertex set $\{1,\dots, 6\}$ whose maximal faces are $\{1,2\}$,  $\{2,3\}$, $\{3,4\}$, $\{4,5\}$, $\{5,1\}$, $\{5,6\}$, $\{2,6\}$.
It has pure height 4. It is a set-theoretic complete intersection, 
since, as we will see, 
\begin{equation}\label{1}I=\sqrt{(x_1x_6,\  x_3x_5,\ x_1x_3+x_2x_4+x_3x_6,\ x_1x_4+x_2x_5+x_4x_6)}.\end{equation}
\noindent
We only need to prove the inclusion $\subset$.
Note that, according to Proposition \ref{main}, 
$$\sqrt{(x_3x_5,\ x_3x_6+x_2x_4, \ x_4x_6+x_2x_5)}=(x_3x_5,\ x_3x_6,\ x_2x_4,\ x_4x_6,\ x_2x_5).$$
Therefore, to prove (\ref{1}) it suffices to show that $x_1x_3$, $x_1x_4$ belong to the right-hand side of (\ref{1}). This is true, because, firstly
\begin{eqnarray*} 
x_1^3x_3^3&=& x_3^2(x_2x_4-x_1x_3)x_1x_6+x_1x_2^2x_3\cdot x_3x_5\\
&&+x_1^2x_3^2(x_1x_3+x_2x_4+x_3x_6)-x_1x_2x_3^2(x_1x_4+x_2x_5+x_4x_6),
\end{eqnarray*}
\noindent
which proves the claim for $x_1x_3$, and secondly,
\begin{eqnarray*} 
x_1^2x_4^2&=&(x_3x_5-x_4^2)x_1x_6+x_1^2\cdot x_3x_5\\
&-&x_1x_5(x_1x_3+x_2x_4+x_3x_6)+x_1x_4(x_1x_4+x_2x_5+x_4x_6).
\end{eqnarray*}
}
\end{example}
Next we apply Proposition \ref{main} to a class of ideals which extends Example \ref{example1}. We will determine the arithmetical rank of these ideals and show that they are not set-theoretic complete intersections except for the ideal studied in that example. 
\begin{example}{\rm Let $m\geq 2$ be an integer, and let $I_m$ be the reduced monomial ideal of $R=K[x_1,\dots, x_{3m+3}]$ generated by the following monomials:
\begin{eqnarray*}
r_1=x_1x_2,\qquad s_n&=&x_{3n-2}x_{3n+2}\\
t_n&=&x_{3n+1}x_{3n+3}\\
u_n&=&x_{3n+1}x_{3n+2}\\
v_n&=&x_{3n-1}x_{3n+3}
\qquad\qquad(n=1,\dots, m).
\end{eqnarray*}
\noindent
We prove that 
\begin{equation}\label{ara}{\rm ara}\,I_m= 2m+1.\end{equation}
\noindent Let $J_m$ be the ideal of $R$ generated by  the following $2m+1$ elements:
$$
x_1x_2,\qquad s_n+t_n,\qquad u_n+v_n\qquad(n=1,\dots, m).
$$
We first show that  
\begin{equation}\label{leq}{\rm ara}\,I_m\leq 2m+1\end{equation}
\noindent
by proving that 
$$I_m=\sqrt{J_m}.$$
\noindent It suffices to show  that $I_m\subset\sqrt{J_m}$, i.e., that $s_n,t_n,u_n, v_n\in \sqrt{J_m}$ for $1\leq n\leq m$. We proceed by finite induction on $n$, $1\leq n\leq m$. For $n=1$ the claim follows by applying Proposition \ref{main} to the following elements:
$$p_0=r_1=x_1x_2,\quad p_1=s_1+t_1=x_1x_5+x_4x_6,\quad p_2=u_1+v_1=x_4x_5+x_2x_6.$$
\noindent
Now suppose that $n>1$ and that the claim is true for $n-1$. Then, in particular, $u_{n-1}=x_{3n-2}x_{3n-1}\in\sqrt{J_m}$. Since $u_{n-1}$ divides $s_nv_n$, $t_n$ divides $u_nv_n$ and $u_n$ divides $s_nt_n$, 
$$p_0=u_{n-1},\quad p_1=s_n+t_n,\quad p_2=u_n+v_n$$
\noindent
fulfil the assumption of Proposition \ref{main}. It follows that $$s_n,t_n,u_n, v_n\in\sqrt{p_0,p_1,p_2}\subset \sqrt{J_m},$$\noindent  which achieves the induction step and proves (\ref{leq}). \newline
We now show the opposite inequality. Let 
$$S=\{x_1\}\cup\{x_{3n+2},\ x_{3n+3},\mid n=1,\dots, m\},$$
\noindent
and set $P=(S)$. Then $P$ is a prime ideal and 
$I_m\subset P$,
since $x_1$ divides $r_1$ and, for all $n=1,\dots, m$, $x_{3n+2}$ divides $s_n$ and $u_n$,  and  $x_{3n+3}$ divides $t_n$ and $v_n$. Ideal $P$ is in fact a minimal prime of $I_m$, because
$$r_1\not\in(S\setminus\{x_1\}),\mbox{ and } u_n\not\in(S\setminus\{x_{3n+2}\}),\ t_n\not\in(S\setminus\{x_{3n+3}\})\quad (n=1,\dots m).$$
\noindent
We have that height\,$P=2m+1$.  By Krull's principal ideal theorem it follows that 
$$2m+1\leq {\rm ara}\,I_m,$$
\noindent
as required. This completes the proof of (\ref{ara}). 
Now let
$$T=\{x_2\}\cup\{x_{3n+1},\ x_{3n+2},\mid n=1,\dots, m-1\}\cup\{x_{3m+1}\},$$
\noindent
and set $Q=(T)$. Then $I_m\subset Q$, since
\begin{list}{}{}
 \item{-} $x_2$ divides $r_1,v_1$;
\item{-} for $n=1,\dots, m-1$, $x_{3n+2}$ divides $s_n$ and $u_n$, and $x_{3n+1}$ divides $t_n$;
\item{-}  for $n=2,\dots, m$, $x_{3(n-1)+2}=x_{3n-1}$ divides $v_n$;
\item{-} $x_{3(m-1)+1}=x_{3m-2}$ divides $s_m$;
\item{-} $x_{3m+1}$ divides $u_m$ and $t_m$.
\end{list}
\noindent
Ideal $Q$ is a minimal prime of $I_m$, because
\begin{list}{}{}
\item{-} $r_1\not\in(T\setminus\{x_2\})$;
\item{-} $t_n\not\in(T\setminus\{x_{3n+1}\})\quad(n=1,\dots, m)$;
\item{-} $v_{n+1}\not\in(T\setminus\{x_{3n+2}\})\quad(n=1,\dots, m-1)$.
\end{list}
\noindent
We have that height\,$Q=1+2(m-1)+1=2m<2m+1$, so that the height of $I_m$ is less than $2m+1$. Hence $I$ is not a set-theoretic complete intersection.   Let us remark that for $m=1$, ideal $I_m$ is, up to a change of variables, the same as ideal $I$ of Example \ref{example1}. }
\end{example}
Finally we present an example of a Stanley-Reisner ideal which is a set-theoretic complete intersection over any field $K$ with more than two elements, but, unlike in all previously examined cases, the minimal set of elements which generate $I$ up to radical strictly depends on $K$. 
\begin{example} {\rm Let $I$ be the Stanley-Reisner ideal of $R=K[x_1, \dots, x_6]$ associated with the Cohen-Macaulay simplicial complex on the vertex set $\{1,\dots, 6\}$ whose maximal faces are $\{1,2\}$,  $\{2,3\}$, $\{3,4\}$, $\{4,5\}$, $\{5,1\}$, $\{4,6\}$, $\{5,6\}$.
 Then 
$$I=(x_1x_3,\ x_1x_4,\ x_1x_6,\ x_2x_4,\ x_2x_5,\ x_2x_6,\ x_3x_5,\ x_3x_6,\ x_4x_5x_6),$$
\noindent
and $I$ is of pure height 4. We show that ara\,$I=4$. First we assume that char\,$K\ne2$. We set 
$$J=(x_1x_4+x_3x_5,\ x_1x_3+x_2x_6+x_4x_5x_6,\ x_1x_6+x_2x_5,\ x_2x_4+x_3x_6)$$
\noindent
 and prove that 
\begin{equation}\label{sqrt}I=\sqrt J.\end{equation}
\noindent
We have that $x_4^2x_5^2x_6^2\in J$, since
\begin{eqnarray*} x_4^2x_5^2x_6^2&=&\frac12x_6(x_6^2-x_1x_4)(x_1x_4+x_3x_5)+x_4x_5x_6(x_1x_3+x_2x_6+x_4x_5x_6)\\
&+&\frac12x_4(x_1x_4-x_6^2)(x_1x_6+x_2x_5)-\frac12x_5(x_1x_4+x_6^2)(x_2x_4+x_3x_6),
\end{eqnarray*}
which implies that $x_4x_5x_6\in\sqrt J$. 
Moreover, $x_3^2x_5^2\in\sqrt J$ (and hence $x_3x_5\in\sqrt J$), because
\begin{eqnarray*}x_3^2x_5^2&=& (x_3x_5-\frac12 x_6^2)(x_1x_4+x_3x_5)-x_4x_5(x_1x_3+x_2x_6)\\
&+&\frac12x_4x_6(x_1x_6+x_2x_5)+\frac12x_5x_6(x_2x_4+x_3x_6),
\end{eqnarray*}
\noindent
where $x_1x_3+x_2x_6\in\sqrt J$.  On the other hand we know that $x_1x_6+x_2x_5\in \sqrt J$. Now Proposition \ref{main} applies to $p_0=x_3x_5$, $p_{11}=x_1x_3$, $p_{12}=x_2x_6$, $p_{21}=x_1x_6$, $p_{22}=x_2x_5$, whence $x_1x_3, x_2x_6, x_1x_6, x_2x_5\in\sqrt J$. Note that $x_3x_5\in \sqrt J$ also implies that $x_1x_4\in \sqrt J$. Finally, from $x_2x_6, x_2x_4+x_3x_6\in\sqrt J$ we deduce, by Lemma \ref{lemma1}, that $x_2x_4, x_3x_6\in\sqrt J$. We have thus proven that $I\subset \sqrt J$; since the opposite inclusion is obvious, (\ref{sqrt}) follows. Now assume that char\,$K=2$, and $K$ has more than two elements. Let $t\in K\setminus\{0,1\}$. This time we set  
$$J=(x_1x_4+x_3x_5,\ x_1x_3+x_2x_6+x_4x_5x_6,\ x_1x_6+tx_2x_5,\ x_2x_4+x_3x_6).$$
\noindent
We show that 
$$I=\sqrt J.$$
\noindent
As above, the claim will follow once we have proven that $x_4x_5x_6, x_3x_5\in\sqrt J$. In fact we have:
\begin{eqnarray*}
x_4^2x_5^2x_6^3&=&\frac1{t+1}\lbrack(x_6^4-tx_1x_2x_3-tx_2^2x_6)(x_1x_4+x_3x_5)+\\
&+&(t+1)x_6(x_4x_5x_6-x_1x_3)(x_1x_3+x_2x_6+x_4x_5x_6)\\
&+&(x_1x_3^2+x_2x_3x_6-x_4x_6^3)(x_1x_6+tx_2x_5)\\
&+&(tx_1^2x_3+tx_1x_2x_6-x_5x_6^3)(x_2x_4+x_3x_6)\rbrack.
\end{eqnarray*}
\noindent
and
\begin{eqnarray*}
x_3^2x_5^2&=&\frac1{t+1}\lbrack((t+1)x_3x_5-x_6^2)(x_1x_4+x_3x_5)-(t+1)x_4x_5(x_1x_3+x_2x_6)\\
&+&x_4x_6(x_1x_6+tx_2x_5)+x_5x_6(x_2x_4+x_3x_6)\rbrack.
\end{eqnarray*}
\noindent
In \cite{B1}, Example 4, we determined a set of four polynomials which generate $I$ up to radical over any algebraically closed field, and, in particular, are independent of the characteristic. They, however, unlike those presented in this example, are not formed by linear combinations of the minimal monomial generators.  
}
\end{example}
\section{The second result}
The ring considered in this section is a polynomial ring $R=K[x_1,\dots, x_N]$, where $K$ is an algebraically closed field.    We will determine the arithmetical rank of certain ideals generated by monomials. We prove a result, based on combinatorial considerations, which generalizes both Lemma \ref{lemma1} and Proposition \ref{main} for this class of ideals and which will allow us to prove the set-theoretic intersection property in various examples. 
\begin{proposition}\label{SVgeneralized} Let $G\subset R$ be a set of monomials. Suppose that there are subsets $S_0\dots, S_r$ of $G$ such that
\begin{list}{}{}
\item[(i)] $\bigcup_{i=0}^rS_i=G$;
\item[(ii)] $S_0$ has exactly one element;
\item[(iii)] the following recursive procedure can always be performed and always comes to an end regardless of the choice of the indeterminate $z$ and the index $j$ at each step.
\begin{list}{}{}
\item[0.] Set $T=S_0$.
\item[1.] Pick an indeterminate $z$ dividing the only element of $T$.
\item[2.] Cancel all monomials divisible by $z$ in every $S_i$. 
\item[3.] If no element of $G$ is left, then end. Else pick an index $j$ such that there is exactly one element left in $S_j$ and set $T=S_j$.
\item[4.] Go to 1.
\end{list}
\end{list}
\noindent
For all $i=0,\dots, r$ we set $q_i=\sum_{\mu\in S_i}\mu$. Then we get
$$\sqrt{(G)}=\sqrt{(q_0,\dots,q_r)}.$$
\end{proposition}
\demo It suffices to show that $\sqrt{(G)}\subset\sqrt{(q_0,\dots,q_r)}.$ We proceed by induction on $r\geq0$. For $r=0$ we have that $(G)=(S_0)=(q_0)$, so that the claim is trivially true. Now suppose that $r>0$ and that the claim is true for all smaller $r$. According to Hilbert's Nullstellensatz, it suffices to show that, whenever all $q_i$ vanish at some ${\bf x}\in K^N$, the same is true for all $\mu\in G$. In the sequel, as long as this does not cause any ambiguity, we will denote a polynomial and its value at ${\bf x}$ by the same symbol. So assume that $q_i=0$ for all $i=0,\dots, r$.  From $q_0=0$ we deduce that one of the indeterminates dividing the only element of $S_0$, say the indeterminate $z$, vanishes. Then all $\mu\in G$ that are divisible by $z$ vanish. Let $\bar G$ be the set of $\mu\in G$ that are not divisible by $z$. We have to show that all $\mu\in \bar G$ vanish. If $\bar G=\emptyset$, then there is nothing to be proven. Otherwise, for all $i=1,\dots, r$, set $\bar S_i=S_i\cap\bar G$.  By assumption we have that $\vert\bar S_j\vert=1$ for some index $j$; up to a change of indices we may assume that $j=1$. Then $\bar G$ and its subsets $\bar S_1,\dots, \bar S_r$ fulfil the assumption of the proposition. For all $i=1,\dots, r$ we set $\bar q_i=\sum_{\mu\in \bar S_i}\mu$. Then by induction $\sqrt{(\bar G)}=\sqrt{(\bar q_1,\dots,\bar q_{r})}.$ Since, by assumption, all $\bar q_i$ vanish, this implies that all $\mu\in\bar G$ vanish, as required. This completes the proof. 
\par\medskip\noindent
\begin{remark}{\rm \begin{list}{}{}
\item[(i)] Note that in the above recursive procedure, the only element left in $S_j$ at step 3 is cancelled as soon as step 2 is performed. We could therefore cancel it right away.
\item[(ii)] In the proof of Proposition \ref{SVgeneralized}, we interpret the recursive procedure  in terms of the vanishing of monomials:  cancelling an indeterminate is the same as supposing that this indeterminate vanishes; cancelling a monomial means concluding that it vanishes. Finishing the procedure means cancelling all monomials of $G$, i.e., concluding that they all vanish. 
\end{list}}
\end{remark}
\begin{remark}{\rm Proposition \ref{SVgeneralized} generalizes Lemma \ref{lemma1} for ideals generated by   monomials over an algebraically closed field. In fact, if $P$ is a set of monomials and and $P_0,\dots, P_r$  are as in the  assumption of Lemma \ref{lemma1}, and we set $G=P$ and  $S_i=P_i$  for all $i=0,\dots, r$, then the assumption of Proposition \ref{SVgeneralized} is fulfilled, as we are going to show next. Let $z$ be any indeterminate dividing the only element of $S_0$. As in the proof of Proposition \ref{SVgeneralized}, for all $i=1,\dots, r$, let $\bar S_i$ denote the set of monomials in $S_i$ which are not divisible by $z$.  Let $j$ be the smallest index $j>0$ such that $\bar S_j$ is not empty. Then the product of each two distinct monomials of $S_j$ is divisible by a monomial of some $S_i$, $i<j$, which is divisible by $z$. Hence all but possibly one of the monomials of $S_j$ are divisible by $z$. Therefore $\vert\bar S_j\vert\leq1$. We conclude by finite descending induction. 
}\end{remark}
\begin{remark}{\rm If the  elements $p_0, p_{11}, p_{12}, p_{21}, p_{22}$ fulfilling the assumption of Proposition \ref{main} are monomials, then it can be easily verified that the assumption of Proposition \ref{SVgeneralized} is fulfilled by $S_0=\{p_0\}$, $S_1=\{p_{11}, p_{12}\}$ and $S_2=\{p_{21}, p_{22}\}$. 
}\end{remark}
\begin{example}\label{i6}{\rm In the ring  $R=K[x_1,\dots, x_6]$ consider the ideal $I$ generated by the following set of squarefree  monomials:
$$ G=\{x_1x_3, x_1x_4, x_1x_5, x_2x_4, x_2x_5, x_2x_6, x_3x_5, x_3x_6,x_4x_6\}.$$
\noindent
We verify that the following subsets of $G$ fulfil the assumption of Proposition \ref{SVgeneralized}:
\begin{eqnarray*} S_0&=&\{x_3x_6\},\\
S_1&=&\{x_1x_4, x_2x_5\},\\
S_2&=&\{x_1x_3, x_2x_4, x_3x_5\},\\
S_3&=&\{x_1x_5, x_2x_6, x_4x_6\}.
\end{eqnarray*}
\noindent 
0. Since $x_3x_6$ is the only element of $S_0$, we cancel it.\newline
1. First pick $x_3$, and cancel $x_1x_3, x_3x_5$. Then $x_2x_4$ is the only element left in $S_2$. We cancel it. 
\newline For the remaining part of the procedure, we have two possible subcases.
\newline
1.1. Pick $x_2$ and cancel $x_2x_5, x_2x_6$. Then $x_1x_4$ is the only element left in $S_1$, and we cancel it. Then we must pick $x_1$ or $x_4$; in the former case we cancel $x_1x_5$, and there is only $x_4x_6$ left in $S_3$;  in the latter case the roles of $x_1x_5$ and $x_4x_6$ are interchanged. Hence, in any case, all monomials are cancelled. 
\newline
1.2. Pick $x_4$ and cancel $x_1x_4, x_4x_6$. Then $x_2x_5$ is the only element left in $S_1$, and we cancel it. Then we must pick $x_2$ or $x_5$; but the first choice takes us back to the previous subcase. In the second case, we cancel $x_1x_5$, and there is only $x_2x_6$ left in $S_3$; we cancel it, and thus all monomials have been cancelled. 
\newline
2. Then pick $x_6$, and cancel $x_2x_6, x_4x_6$. Then $x_1x_5$ is the only element left in $S_3$. We cancel it. 
\newline For the sequel of the procedure, we again have two possible subcases.
\newline
2.1. Pick $x_1$ and cancel $x_1x_3, x_1x_4$. Then $x_2x_5$ is the only element left in $S_1$, and we cancel it. Then we must pick $x_2$ or $x_5$; in the former case we cancel $x_2x_4$, and there is only $x_3x_5$ left in $S_2$;  in the latter case the roles of $x_2x_4$ and $x_3x_5$ are interchanged. Hence, in any case, all monomials are cancelled. 
\newline
2.2. Pick $x_5$ and cancel $x_2x_5, x_3x_5$. Then $x_1x_4$ is the only element left in $S_1$, and we cancel it. Then we must pick $x_1$ or $x_4$; but the first choice takes us back to the previous subcase. In the second case, we cancel $x_2x_4$, and there is only $x_1x_3$ left in $S_2$; we cancel it, and thus all monomials have been cancelled. 
\newline
Therefore 
\begin{equation}\label{terai}I=\sqrt{(x_3x_6,\ x_1x_4+x_2x_5,\ x_1x_3+x_2x_4+x_3x_5,\ x_1x_5+x_2x_6+x_4x_6)}.\end{equation}
\noindent
This equality was first discovered by Naoki Terai by other means and communicated to the author. Four polynomials generating $I$ up to radical were already given in \cite{B1}, Example 2, but they are of a much more complicated form. The attempt to place (\ref{terai}) in a more general framework was the motivation for Proposition \ref{SVgeneralized}. The ideal $I$  is the Stanley-Reisner ideal of the  (Gorenstein) simplicial complex on the vertex set $\{1,\dots, 6\}$ whose maximal faces are $\{1, 2\}$, $\{2, 3\}$, $\{3, 4\}$, $\{4, 5\}$, $\{5, 6\}$, $\{6, 1\}$.  The ideal $I$ has pure height 4. Hence $I$ is a set-theoretic complete intersection.    
}\end{example}
\begin{remark}{\rm The ideals of Example \ref{example1} and Example \ref{i6} are Gorenstein. According to \cite{BH}, Corollary 5.5.6, more generally,  the Stanley-Reisner ideal of the simplicial complex on $N$ vertices whose maximal faces are the edges of a simple closed $N$-gon is always Gorenstein. For $N=4$ this ideal is generated by $x_1x_3$ and $x_2x_4$ and is therefore a complete intersection.  We are not able to say whether the Stanley-Reisner ideal associated with an $N$-gon is a set-theoretic complete intersection for $N\geq7$.
 }
\end{remark}
In addition to Example \ref{i6} we present two more examples  of simplicial complexes on $N=6$ vertices whose maximal faces have cardinality $d=2$ and whose Stanley-Reisner ideals are set-theoretic complete intersections. The combinatorial verifications are left to the reader.  
\begin{example} {\rm Let $R=K[x_1,\dots,  x_6]$. The Stanley-Reisner ideal of the Cohen-Macaulay simplicial complex on the vertex set $\{1,\dots, 6\}$ whose maximal faces are $\{1,2\}$,  $\{2,3\}$, $\{3,4\}$, $\{4,5\}$, $\{5,6\}$, $\{6,1\}$, $\{2,5\}$ is 

$$I=(x_1x_3,\ x_1x_4,\ x_1x_5,\ x_2x_4,\ x_2x_6,\ x_3x_5,\ x_3x_6,\ x_4x_6).$$ 
\noindent
From Proposition \ref{SVgeneralized} we can deduce that
$$I=\sqrt{(x_1x_4,\ x_3x_6,\ x_1x_3+x_1x_5+x_2x_4,\ x_2x_6+x_3x_5+x_4x_6)}.$$ 
}\end{example}
\begin{example}{\rm 
The Stanley-Reisner ideal of the Cohen-Macaulay simplicial complex  on the vertex set $\{1,\dots, 6\}$ whose maximal faces are $\{1,2\}$,  $\{2,3\}$, $\{3,4\}$, $\{4,5\}$, $\{5,6\}$, $\{6,1\}$, $\{2,5\}$, $\{2,6\}$, $\{1,5\}$
is
$$I=(x_1x_3,\ x_1x_4,\ x_2x_4,\ x_3x_5,\ x_3x_6,\ x_4x_6,\ x_1x_2x_5,\ x_1x_2x_6,\ x_1x_5x_6,\ x_2x_5x_6).$$ 
\noindent
From Proposition \ref{SVgeneralized} we can deduce that
\begin{eqnarray*}I&=&\sqrt{(x_1x_4, x_3x_6, x_1x_3+x_2x_4+x_1x_2x_5},\\
 &&\qquad\qquad\qquad\overline{x_3x_5+x_4x_6+x_1x_2x_6+x_1x_5x_6+x_2x_5x_6)}
\end{eqnarray*}
}
\end{example}
Many more examples could be provided. We only give one for $N=8$, $d=2$.
\begin{example} {\rm Let $R=K[x_1,\dots,  x_8]$. The Stanley-Reisner ideal of the Cohen-Macaulay simplicial complex  on the vertex set $\{1,\dots, 8\}$ whose maximal faces are $\{1,2\}$,  $\{2,3\}$, $\{3,4\}$, $\{4,5\}$, $\{5,1\}$,  $\{5,6\}$, $\{6,7\}$, $\{7,8\}$, $\{8,4\}$ is
\begin{eqnarray*}I&=&(x_1x_3,\ x_1x_4,\ x_1x_6,\ x_1x_7,\ x_1x_8,\\
&&\hphantom{(}x_2x_4,\ x_2x_5,\ x_2x_6,\ x_2x_7,\ x_2x_8,\\
&&\hphantom{(}x_3x_5,\ x_3x_6,\ x_3x_7,\ x_3x_8,\ x_4x_6,\ x_4x_7,\ x_5x_7,\ x_5x_8,\ x_6x_8).
\end{eqnarray*} 
\noindent
From Proposition \ref{SVgeneralized} we can deduce that
\begin{eqnarray*}
I&=&\sqrt{(x_3x_6,\ x_1x_8+x_2x_7,\ x_1x_3+x_2x_8+x_3x_7,\ x_1x_7+x_2x_6+x_6x_8},\\
&&\overline{x_1x_4+x_2x_5+x_3x_8+x_4x_7+x_5x_8,\  x_1x_6+x_2x_4+x_3x_5+x_4x_6+x_5x_7)}.
\end{eqnarray*}
Since $I$ is an ideal of pure height 6, we conclude that it is a set-theoretic complete intersection.
}\end{example}

\vskip.03truecm
\begin{center}{\bf ACKNOWLEDGEMENTS}
\end{center}
\noindent
The author is indebted to Naoki Terai for the interesting example that motivated one of the main results of this paper. 
\end{document}